\documentclass[12pt,a4paper]{amsart}

\title[Tree Product Lattices]{Tiling systems and homology of lattices in tree products}
\author{Guyan Robertson}
\address{School of Mathematics and Statistics, University of Newcastle, NE1 7RU, U.K.}
\email{a.g.robertson@newcastle.ac.uk}

\usepackage{amsmath}
\usepackage{amscd}
\usepackage{amssymb}

\input{prepictex}
\input{pictex}
\input{postpictex}
\chardef\bslash=`\\ 

\makeatletter
\def\verbatim{\interlinepenalty\@M \@verbatim
  \leftskip\@totalleftmargin\advance\leftskip2pc
  \frenchspacing\@vobeyspaces \@xverbatim}
\makeatother
\newtheorem{theorem}{Theorem}[section]

\newtheorem{lemma}[theorem]{Lemma}

\theoremstyle{definition}

\newcommand{\ovl}{\overline}
\newcommand{\n}{\noindent}
\newcommand{\cl}[1]{{\mathcal{#1}}}
\newcommand{\bb}[1]{{\mathbb{#1}}}
\newcommand{\fk}[1]{{\mathfrak{#1}}}

\newcounter{picture}

\DeclareMathOperator{\rank}{rank}

\DeclareMathOperator{\range}{range}

\newcommand{\e}{{\varepsilon}}

\newcommand{\vf}{{\varphi}}
\newcommand{\aut}{{\text{\rm Aut}}}

\newcommand{\pgl}{{\text{\rm{PGL}}}}

\newcommand{\0}{{\bf 0}}
\newcommand{\1}{{\bf 1}}

\begin{document}

\begin{abstract}
Let $\Gamma$ be a torsion free cocompact lattice in $\aut(\cl T_1)\times\aut(\cl T_2)$,
where $\cl T_1$, $\cl T_2$ are trees whose vertices all have degree at least three.
The group $H_2(\Gamma, \bb Z)$ is determined explicitly in terms of an associated
2-dimensional tiling system. 
It follows that under appropriate conditions the crossed product $C^*$-algebra $\cl A$ associated with the action of $\Gamma$ on the boundary of $\cl T_1\times \cl T_2$ satisfies $\rank K_0(\cl A) = 2\cdot\rank H_2(\Gamma, \bb Z)$.
\end{abstract}

\maketitle

\section{Introduction}

This article is motivated by the problem of calculating the K-theory of certain crossed product
$C^*$-algebras $\cl A(\Gamma,\partial\Delta)$, where $\Gamma$ is a higher rank lattice
acting on an affine building $\Delta$ with boundary $\partial\Delta$. 
Here we examine the case where $\Delta$ is a product of trees. We determine the K-theory rationally, thereby proving some conjectures in \cite{kr}.

Let $\cl T_1$ and $\cl T_2$ be locally finite trees whose vertices all have degree at least three.
Consider the direct product $\Delta = \cl T_1\times \cl T_2$ as a two dimensional cell complex.
Let $\Gamma$ be a discrete subgroup of $\aut(\cl T_1)\times\aut(\cl T_2)$ which acts freely and cocompactly on $\Delta$.
Associated with the action $(\Gamma, \Delta)$ is a tiling system whose set of tiles is the set $\fk R$ of ``directed'' 2-cells of $\Gamma\backslash\Delta$. There are vertical and horizontal adjacency rules $tHs$ and $tVs$ between tiles $t, s \in\fk R$ illustrated below. 
Precise definitions will be given in Section \ref{prod}.
\refstepcounter{picture}
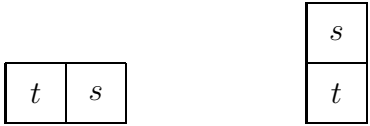
\begin{figure}[htbp]\label{adjacency}
\hfil
\centerline{
\beginpicture
\setcoordinatesystem units <.8cm, .8cm>  
\setplotarea  x from -6 to 6,  y from 0 to 1.5
\putrule from -3 0 to -1 0
\putrule from -3 1 to -1 1
\putrule from -3 0 to -3 1
\putrule from -2 0 to  -2 1
\putrule from -1 0 to  -1 1
\putrule from 2 0 to  3 0
\putrule from 2 1 to  3 1
\putrule from 2 2 to  3 2
\putrule from 2 0 to  2 2
\putrule from 3 0 to  3 2
\put {$t$} []     at   -2.5 0.5
\put {$s$} at   -1.5 0.5
\put {$t$} []     at   2.5 0.5
\put {$s$} at   2.5  1.5
\endpicture
}
\hfil
\end{figure}

\noindent
There are homomorphisms $T_1, T_2 : \bb Z \fk R \to \bb Z \fk R$ defined by
\begin{equation*}
  T_1t = \sum_{tHs}s, \qquad
  T_2t = \sum_{tVs}s\,.
\end{equation*}
Consider the homomorphism $\bb Z \fk R \to \bb Z \fk R  \oplus \bb Z \fk R$ given by
\begin{equation*}
\left(\begin{matrix}T_1-I\\T_2-I\end{matrix}\right) :
 t\mapsto (T_1t-t)\oplus (T_2t-t).  
 \end{equation*}
The main result of this article it the following Theorem, which is formulated more precisely in Theorem \ref{main}.
\begin{theorem}\label{main0} There is an isomorphism 
\begin{equation}\label{p}
H_2(\Gamma, \bb Z)\cong\ker \left(\begin{smallmatrix}T_1-I\\T_2-I\end{smallmatrix}\right).
\end{equation}
\end{theorem}
The proof of (\ref{p}) is elementary, but care is needed because the right hand side is defined in terms of ``directed'' 2-cells rather than geometric 2-cells.
A square complex $X$ is VH-T if every vertex link is a complete bipartite graph 
and if there is a partition of the set of edges into vertical and horizontal, which agrees
with the bipartition of the graph on every link \cite{bm}. The universal covering space $\Delta$ of a VH-T complex $X$ is a product of trees $\cl T_1\times \cl T_2$ and the fundamental group $\Gamma$ of $X$ is a subgroup of $\aut(\cl T_1)\times\aut(\cl T_2)$
which acts freely and cocompactly on $\cl T_1\times \cl T_2$. Conversely, every finite VH-T complex arises
in this way from a free cocompact action of a group $\Gamma$ on a product of trees. Recall that a discrete group which acts freely on a CAT(0) space is necessarily torsion free. 

The group $\Gamma$ acts on the (maximal) boundary $\partial\Delta$ of $\Delta$, which is 
the set of chambers of the spherical building at infinity, endowed with an appropriate topology \cite{kr}.
This boundary may be identified with a direct product of Gromov boundaries $\partial\cl T_1\times \partial\cl T_2$.
The boundary action $(\Gamma,\partial\Delta)$ gives rise to a crossed product $C^*$-algebra 
$\cl A(\Gamma,\partial\Delta)=C_{\bb C}(\partial\Delta) \rtimes \Gamma$ as described in \cite{kr}.

If $p$ is prime then $\pgl_2(\bb Q_p)$ acts on its Bruhat-Tits tree $\cl T_{p+1}$, which is a 
homogeneous tree of degree $p+1$.
If $p, \ell$ are prime then the group $\pgl_2(\bb Q_p)\times\pgl_2(\bb Q_\ell)$ acts on the $\Delta=\cl T_{p+1} \times \cl T_{\ell+1}$.
Let $\Gamma$ be a torsion free irreducible lattice in $\pgl_2(\bb Q_p)\times\pgl_2(\bb Q_\ell)$.
Then $\cl A(\Gamma,\partial\Delta)$ is a higher rank
Cuntz-Krieger algebra and fits into the general theory developed in \cite{rs1,rs2}.
In particular, it is classified up to isomorphism by its K-theory. It is a consequence of Theorem \ref{main0}
(see Section \ref{five}) that
\begin{equation}\label{principal}
\rank K_0(\cl A(\Gamma,\partial\Delta))= 2\cdot\rank H_2(\Gamma, \bb Z).
\end{equation} 
This proves a conjecture in \cite{kr}. 
The normal subgroup theorem \cite[IV, Theorem (4.9)]{mar} implies that $H_1(\Gamma, \bb Z)$ is a finite group.
Equation (\ref{principal}) can therefore be expressed as
\begin{equation*}
\chi(\Gamma) =1+ \frac{1}{2}\rank K_0(\cl A(\Gamma,\partial\Delta)).
\end{equation*}
One easily calculates that $\chi(\Gamma)=\frac{(p-1)(\ell-1)}{4}|X^0|$, where $|X^0|$
is the number of vertices of $X$. Therefore the rank of $K_0(\cl A(\Gamma,\partial\Delta))$ can be expressed explicitly in terms of $p,\ell$ and $|X^0|$. Examples are constructed in \cite[Section 3]{moz}, where $p,\ell \equiv 1$ (mod 4) are two distinct primes. 

\section{Products of trees and their automorphisms.}\label{prod}

If $\cl T$ is a tree, there is a type map $\tau$ defined on the
vertex set of $\cl T$, taking values in $\bb Z/2\bb Z$.
Two vertices have the same type if and only if the distance between them is even. Any automorphism $g$ of $\cl T$ preserves distances between vertices, and so there exists $i\in\bb Z/2\bb Z$ (depending on $g$) such that $\tau(gv)= \tau(v)+i$, for every vertex $v$.

Suppose that $\Delta$ is the $2$ dimensional cell complex associated with  a product $\cl T_1\times \cl T_2$ of trees.
Let $\Delta^k$ denote the set of $k$-cells in $\Delta$ for $k=0,1,2$. The 0-cells are vertices and the 2-cells are geometric squares.
Denote by $u=(u_1,u_2)$ a generic vertex of $\Delta$.
There is a type map~$\tau$ on $\Delta^0$ defined by
\begin{equation*}
\tau(u)=
(\tau(u_1),\tau(u_2))\in\bb Z/2\bb Z\times\bb Z/2\bb Z.
\end{equation*}
Any 2-cell $\delta \in \Delta^2$ has one vertex of each type.
For every $g\in \aut \cl T_1 \times \aut \cl T_2$ there exists $(k,l)\in\bb Z/2\bb Z\times\bb Z/2\bb Z$
such that, for each vertex $u$,
\begin{equation}\label{type}
\tau(gu)=(\tau(u_1)+k,\tau(u_2)+l).
\end{equation}   
Let $\Gamma < \aut \cl T_1 \times \aut \cl T_2$ be a torsion free discrete group acting cocompactly on $\Delta$. Then $X=\Gamma\backslash\Delta$ is a finite cell complex with universal covering $\Delta$. Let $X^k$ denote the set of $k$-cells of $X$ for $k=0,1,2$.

\begin{figure}[htbp]
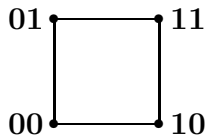
\label{modelsquare}
\hfil
\centerline{
\beginpicture
\setcoordinatesystem units <0.7cm, 0.7cm>  
\setplotarea  x from -6 to 6,  y from -1 to 1.2
\putrule from -1 -1 to 1 -1
\putrule from -1  1 to 1  1
\putrule from -1 -1 to -1 1
\putrule from  1  -1 to 1  1
\put {$_{\bullet}$} at 1 1
\put {$_{\bullet}$} at 1 -1
\put {$_{\bullet}$} at -1 1
\put {$_{\bullet}$} at -1 -1
\put {${\1\0}$}[l] at 1.2 -1
\put {${\0\1}$}[r] at -1.2 1
\put {${\0\0}$}[r] at -1.2 -1
\put {${\1\1}$}[l] at 1.2  1
\endpicture
}
\hfil
\caption{The model square $\sigma$.}
\end{figure}

The first step is to formalize the notion of a directed square in $X$. We modify the terminology of \cite[Section 1]{bm}, in order to fit with \cite{rs1,rs2,kr}. Let $\sigma$ be a model typed square with vertices $\0\0, \0\1, \1\0, \1\1$, as illustrated in Figure \ref{modelsquare}. Assume that the vertex $\bf ij$ of $\sigma$ has type $$\tau({\bf ij})= (i,j) \in \bb Z/2\bb Z\times\bb Z/2\bb Z.$$ 

\refstepcounter{picture}
The vertical and horizontal reflections $v,h$ of $\sigma$ are the involutions satisfying
$v({\bf 00}) = {\bf 01}, v({\bf 10}) = {\bf 11}, h({\bf 00}) = {\bf 10}, h({\bf 01}) = {\bf 11}$.
An isometry $r: \sigma \to \Delta$ is said to be \emph{type rotating} if there exists
$(k,l)\in\bb Z/2\bb Z\times\bb Z/2\bb Z$
such that, for each vertex $\bf ij$ of $\sigma$
$$
\tau(r({\bf ij}))=(i+k, j+l).
$$
Let  $R$ denote the set of type rotating isometries $r: \sigma \to \Delta$.
If $g\in \aut \cl T_1 \times \aut \cl T_2$ and $r\in R$ then it follows from (\ref{type}) that $g\circ r \in R$.
If $\delta^2 \in \Delta^2$ then for each $(k,l)\in\bb Z/2\bb Z\times\bb Z/2\bb Z$ there is a unique $r \in R$ such that $r(\sigma)=\delta^2$ and $r(\0\0)$ has type $(k,l)$. Therefore each geometric square $\delta^2 \in \Delta^2$ is the image of each of the four elements of $\{ r\in R\ ; \ r(\sigma)=\delta^2\}$ under the map $r\mapsto r(\sigma)$.
The next lemma records this observation.

\begin{lemma}\label{nto1}
The map $r\mapsto r(\sigma)$ from $R$ to $\Delta^2$ is $4$-to-$1$.
\end{lemma}

Let $\fk R = \Gamma\backslash R$ and call $\fk R$  the set of {\em directed squares} of $X=\Gamma\backslash\Delta$.
There is a commutative diagram
\begin{equation*}
\begin{CD}
R                     @>r\mapsto r(\sigma)>>   \Delta^2\\
@VVV                        @VVV \\
\fk R                  @>\eta>>      X^2
\end{CD}
\end{equation*}
where the vertical arrows represent quotient maps and  $\eta$ is defined by $\eta(\Gamma r)=\Gamma.r(\sigma)$.
The next result makes precise the fact that each geometric square in $X^2$ corresponds to exactly four directed squares.

\begin{lemma}\label{epsilon}
The map $\eta : \fk R \to X^2$ is surjective and $4$-to-$1$.
\end{lemma}

\begin{proof}
Fix $\delta^2\in R$.
By Lemma \ref{nto1}, the set $\{ r\in R\ ; \ r(\sigma)=\delta^2\}=\{r_1, r_2, r_3, r_4\}$ contains precisely $4$ elements.
Since $\Gamma$ acts freely on $\Delta$, the set $\{\Gamma r_1, \Gamma r_2, \Gamma r_3, \Gamma r_4\} \subset \fk R$ also contains precisely four elements, each of which maps to $\Gamma\delta^2$ under $\eta$.
Now suppose that $\eta(\Gamma r)=\Gamma\delta^2$ for some $r\in R$.
Then $\gamma r(\sigma)=\delta^2$ for some $\gamma\in\Gamma$. Thus $\gamma r\in \{r_1, r_2, r_3, r_4\}$ and
$\Gamma r \in \{\Gamma r_1, \Gamma r_2, \Gamma r_3, \Gamma r_4\}$.
 This proves that $\eta$ is $4$-to-$1$.
\end{proof}

The vertical and horizontal reflections $v, h$ of the model square $\sigma$ act on $\fk R$ and generate a group
$\Sigma\cong \bb Z/2\bb Z\times\bb Z/2\bb Z$ of symmetries of $\fk R$. The $\Sigma$-orbit of each $r\in \fk R$ contains four elements. Choose once and for all a subset 
$\fk R^+ \subset \fk R$ containing precisely one element from each $\Sigma$-orbit.
The map $\eta$ restricts to a 1-1 correspondence between $\fk R^+$ and the set of geometric squares $X^2$.
For each $\phi\in \Sigma-\{1\}$, let $\fk R^\phi$ denote the image of  $\fk R^+$ under $\phi$. Then $\fk R$ may be expressed as a disjoint union
\begin{equation*}
\fk R = \fk R^+\cup  \fk R^v\cup  \fk R^h\cup  \fk R^{vh}.
\end{equation*}
Now we formalize the notion of horizontal and vertical directed edges in $X$.
Consider the two directed edges $[\0\0, \1\0], [\0\0, \0\1]$ of the model square $\sigma$.

\refstepcounter{picture}
\begin{figure}[htbp]
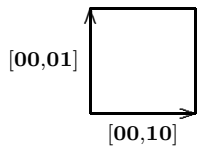
\label{modeledges}
\hfil
\centerline{
\beginpicture
\setcoordinatesystem units <0.7cm, 0.7cm>  
\setplotarea  x from -6 to 6,  y from -1 to 1.2
\putrule from -1 -1 to 1 -1
\putrule from -1  1 to 1  1
\putrule from -1 -1 to -1 1
\putrule from  1  -1 to 1  1
\put {$_{[\0\0, \0\1]}$} [r]     at   -1.2 0
\put {$_{[\0\0, \1\0]}$} [t]     at   0 -1.2
\arrow <6pt> [.3,.67] from    0.8 -1   to    1 -1
\arrow <6pt> [.3,.67] from    -1  0.8   to  -1 1 
\endpicture
}
\hfil
\caption{Directed edges of the model square $\sigma$.}
\end{figure}
Let  $A$ be the set of type rotating isometries $r: [\0\0,\1\0] \to \Delta$, and
let  $B$ be the set of type rotating isometries $r: [\0\0,\0\1] \to \Delta$.
There is a natural 2-to-1 mapping $r\mapsto \range r$, from $A \cup B$ onto $\Delta^1$.
Let $\fk A = \Gamma\backslash A$ and $\fk B = \Gamma\backslash B$.
Call $\fk A, \fk B$ the sets of horizontal and vertical \textit{directed edges} of $X=\Gamma\backslash\Delta$. Let $\cl E = \fk A \cup \fk B$, the set of all directed edges of $X$. 

If $a=\Gamma r\in \fk A$, let $o(a)=\Gamma r(\0\0)\in X^0$ and $t(a)= \Gamma r(\1\0)\in X^0$, the \textit{origin} and \textit{terminus} of the directed edge $a$. Similarly, if $b=\Gamma r\in \fk B$, let $o(b)=\Gamma r(\0\0)\in X^0$ and $t(b)= \Gamma r(\0\1)\in X^0$. Note that it is possible that $o(e)=t(e)$.

A straightforward analogue of Lemma \ref{epsilon} shows that each geometric edge in $X^1$ is the image of each of two directed edges. 
The horizontal and vertical reflections on $\sigma$ induce an inversion on $\cl E$,
denoted by $e\mapsto \ovl e$, with the property that $\ovl{\ovl e} = e$ and $o(e)=t(\ovl e)$.
The pair $(\cl E, X^0)$ is thus a graph in the sense of \cite{ser2}.
Choose once and for all an orientation of this graph: that is a subset
$\cl E^+$ of $\cl E$, with $\cl E = \cl E^+ \sqcup \overline{\cl E^+}$.
Write $\fk A^+ = \fk A \cap \cl E^+$ and $\fk B^+ = \fk B \cap \cl E^+$. The images of $\fk A$ [respectively $\fk B$] in $X^1$
are the edges the \emph{horizontal [vertical] 1-skeleton} $X^1_h$ [$X^1_v$].

\begin{lemma}
There is a well defined injective map 
$$t\mapsto (a(t), b(t)) : \fk R \to \fk A \times \fk B$$
which is surjective if $X$ has one vertex.
\end{lemma}

\refstepcounter{picture}
\begin{figure}[htbp]
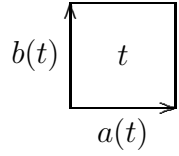
\label{directededges}
\hfil
\centerline{
\beginpicture
\setcoordinatesystem units <0.7cm, 0.7cm>  
\setplotarea  x from -6 to 6,  y from -1 to 1.2
\putrule from -1 -1 to 1 -1
\putrule from -1  1 to 1  1
\putrule from -1 -1 to -1 1
\putrule from  1  -1 to 1  1
\put {$b(t)$} [r]     at   -1.2 0
\put {$a(t)$} [t]     at   0 -1.2
\put {$t$}     at   0  0
\arrow <6pt> [.3,.67] from    0.8 -1   to    1 -1
\arrow <6pt> [.3,.67] from    -1  0.8   to  -1 1 
\endpicture
}
\hfil
\caption{Directed edges in $X$.}
\end{figure}

\begin{proof}
The map $r\mapsto (r|_{[\0\0, \1\0]}, r|_{[\0\0, \0\1]}) : R \to A \times B$ is injective because
each geometric square of $\Delta$ is uniquely determined by any two edges containing a common vertex.

If $t=\Gamma r\in \fk R$ then define
$$a(t)= \Gamma r|_{[\0\0, \1\0]},\quad b(t)= \Gamma r|_{[\0\0, \0\1]}.$$
Using the fact that $\Gamma$ acts freely on $\Delta$ it is easy to see that the map
$t\mapsto (a(t), b(t))$ is injective. 

If $X$ has one vertex, then any two elements $a\in \fk A$, $b\in \fk B$ are represented by type
rotating isometries $r_1: [\0\0,\1\0] \to \Delta$, $r_2: [\0\0,\0\1] \to \Delta$ with 
$r_1(\0\0)=r_2(\0\0)$. The isometries $r_1, r_2$ are restrictions of an isometry $r\in R$, which defines 
an element $t=\Gamma r \in \fk R$ with $a=a(t)$ and $b=b(t)$.
\end{proof}

If $t=\Gamma r\in \fk R$, define directed edges $a'(t)\in \fk A, b'(t)\in \fk B$ opposite to $a(t), b(t)$,
as follows.
\begin{align*}
  a'(t) &= \Gamma(r\circ v|_{[\0\0, \1\0]}), \\
  b'(t) &= \Gamma(r\circ h|_{[\0\0, \0\1]}). 
\end{align*}
\refstepcounter{picture}
\begin{figure}[htbp]
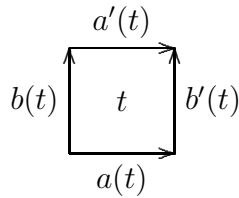
\label{oppedges}
\hfil
\centerline{
\beginpicture
\setcoordinatesystem units <0.7cm, 0.7cm>  
\setplotarea  x from -6 to 6,  y from -1 to 1.2
\putrule from -1 -1 to 1 -1
\putrule from -1  1 to 1  1
\putrule from -1 -1 to -1 1
\putrule from  1  -1 to 1  1
\put {$b(t)$} [r]     at   -1.2 0
\put {$b'(t)$} [l]     at    1.2 0
\put {$a(t)$} [t]     at   0 -1.2
\put {$a'(t)$} [b]     at   0  1.2
\put {$t$}     at   0  0
\arrow <6pt> [.3,.67] from    0.8 -1   to    1 -1
\arrow <6pt> [.3,.67] from    0.8  1   to  1  1
\arrow <6pt> [.3,.67] from    -1  0.8   to  -1 1 
\arrow <6pt> [.3,.67] from   1  0.8    to   1  1
\endpicture
}
\hfil
\caption{Opposite edges.}
\end{figure}
In other words
\begin{equation}\label{flip}
  a'(t)=a(t^v); \quad b'(t)=b(t^h).
\end{equation} 

\section{Some related graphs}

Associated to the VH-T complex $X$ are two graphs whose vertices are directed edges of $X$. Denote
by $\cl G_v(\fk A)$ the graph whose vertex set is $\fk A$ and whose edge set is $\fk R$,
with origin and terminus maps defined by $t\mapsto a(t)$ and $t\mapsto a'(t)$ respectively. 
Similarly $\cl G_h(\fk B)$ is the graph whose vertex set is $\fk B$ and whose edge set is $\fk R$,
with the origin and terminus maps defined by $t\mapsto b(t)$ and $t\mapsto b'(t)$. 
\refstepcounter{picture}
\begin{figure}[htbp]\label{lineedges}
\hfil
\centerline{
\beginpicture
\setcoordinatesystem units <0.7cm, 0.7cm>  
\setplotarea  x from -2 to 8,  y from -1 to 1.2
\putrule from -1 -1 to 1 -1
\putrule from -1  1 to 1  1
\putrule from -1 -1 to -1 1
\putrule from  1  -1 to 1  1
\put {$a'(t)$} [b]     at  0  1.2
\put {$a(t)$} [t]     at   0 -1.2
\put {$t$}     at   0  0
\arrow <6pt> [.3,.67] from    0.8 -1   to    1 -1
\arrow <6pt> [.3,.67] from   0.8 1   to    1 1
\putrule from 5 -1 to 7 -1
\putrule from 5  1 to 7  1
\putrule from 5 -1 to 5 1
\putrule from  7  -1 to 7  1
\put {$b(t)$} [r]     at   4.8  0
\put {$b'(t)$} [l]     at  7.2  0
\put {$t$}     at   6  0
\arrow <6pt> [.3,.67] from    7  0.8   to  7 1 
\arrow <6pt> [.3,.67] from    5  0.8   to  5 1 
\endpicture
}
\hfil
\caption{Edges of $\cl G_v(\fk A)$ and $\cl G_h(\fk B)$.}
\end{figure}

Now define two \textit{directed} graphs whose vertices are elements of $\fk R$.
The ``horizontal'' graph $\cl G_h(\fk R)$ has vertex set $\fk R$. A directed edge
$[t,s]$ is defined as follows. 
Consider the model rectangle $H$ made up of two adjacent squares with vertices $\{(i, j) \in \bb Z^2 : i=0,1,2, j=0,1\}$ where the vertex $(i,j)$ has type $(i+2\bb Z, j+ 2\bb Z)$.
The model square $\sigma$ of Figure \ref{modelsquare} is considered as the left hand square of $H$. 
\refstepcounter{picture}
\begin{figure}[htbp]
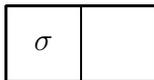
\label{trans1}
\hfil
\centerline{
\beginpicture
\setcoordinatesystem units <1cm, 1cm>  
\setplotarea  x from -1.5 to 1.5,  y from -0.2 to 1.2
\putrule from -1 0 to 1 0
\putrule from -1 1 to 1 1
\putrule from -1 0 to -1 1
\putrule from 0 0 to  0 1
\putrule from 1 0 to  1 1
\put {$\sigma$} at   -0.5  0.5
\endpicture
}
\hfil
\caption{The model rectangle $H$.}
\end{figure}

An isometry $r: H \to \Delta$ is said to be type rotating if there exists
$(k,l)\in\bb Z/2\bb Z\times\bb Z/2\bb Z$
such that, for each vertex $(i,j)$ of $H$,
$\tau(r((i,j)))=(i+k, j+l)$.
A \textit{directed edge} of $\cl G_h(\fk R)$ is $\Gamma r$ where $r : H \to \Delta$ be a type rotating isometry. The \textit{origin} of $\Gamma r$ is $t= \Gamma r_1$, where
$r_1=r|_\sigma$ and the \textit{terminus} of $\Gamma r$ is $s=\Gamma r_2$, where $r_2 : \sigma \to \Delta$ is
defined by $r_2(i,j)=r(i+1, j)$. 
There is a similar definition for the ``vertical'' graph $\cl G_v(\fk R)$ with vertex set $\fk R$.
Edges $[t,s]$ of $\cl G_h(\fk R)$ and $\cl G_v(\fk R)$ are illustrated in Figure \ref{ate}, by the ranges of representative isometries.
\refstepcounter{picture}\label{ate}
\begin{figure}[htbp]\label{HM}
\hfil
\centerline{
\beginpicture
\setcoordinatesystem units <1cm, 1cm>  
\setplotarea  x from -6 to 6,  y from -1 to 2.2
\putrule from -3 0 to -1 0
\putrule from -3 1 to -1 1
\putrule from -3 0 to -3 1
\putrule from -2 0 to  -2 1
\putrule from -1 0 to  -1 1
\putrule from 2 0 to  3 0
\putrule from 2 1 to  3 1
\putrule from 2 2 to  3 2
\putrule from 2 0 to  2 2
\putrule from 3 0 to  3 2
\put {An edge of $\cl G_h(\fk R)$} []     at   -2 -1
\put {An edge of $\cl G_v(\fk R)$} at  2.5 -1
\put {$t$} []     at   -2.5 0.5
\put {$s$} at   -1.5 0.5
\put {$t$} []     at   2.5 0.5
\put {$s$} at   2.5  1.5
\endpicture
}
\hfil
\caption{}
\end{figure}

Since $\Gamma$ acts freely on $\Delta$, it is easy to see that the existence of a directed 
edge  $[t,s]$ of $\cl G_h(\fk R)$ with origin $t\in\fk R$ and terminus $s\in \fk R$ is equivalent to
\begin{equation}\label{e2}
  b(s)=b'(t),\quad  s\ne t^h.
\end{equation}
Similarly the existence of a directed 
edge  $[t,s]$ of $\cl G_v(\fk R)$, with origin $t\in\fk R$ and terminus $s\in \fk R$ is equivalent to
\begin{equation}\label{e3}
  a(s)=a'(t),\quad  s\ne t^v.
\end{equation}
The next Lemma will be used later.
Recall that a lattice $\Gamma$ in $\pgl_2(\bb Q_p)\times\pgl_2(\bb Q_\ell)$ is automatically cocompact \cite[IX Proposition 3.7)]{mar}.
\begin{lemma}\label{nasc}
  If $p, \ell$ are prime and $\Gamma$ is a torsion free irreducible lattice in $\pgl_2(\bb Q_p)\times\pgl_2(\bb Q_\ell)$
  acting on the corresponding product of trees, then the directed graphs  $\cl G_h(\fk R)$, $\cl G_v(\fk R)$ are connected.
\end{lemma}

\begin{proof}
  This follows from \cite[Proposition 2.15]{moz}, using the topological transitivity of an associated shift system. The proof uses the Howe-Moore theorem for $p$-adic semisimple groups and is explained in \cite[Lemma 2]{moz2}. 
\end{proof}

\section{Tilings and $H_2(\Gamma, \bb Z)$}

Throughout this section, $\cl T_1$ and $\cl T_2$ are locally finite trees whose vertices all have degree at least three.
The group $\Gamma$ acts freely and cocompactly on the $2$ dimensional cell complex $\Delta = \cl T_1\times \cl T_2$
and we continue to use the notation introduced in the preceding sections.

For $t, s \in \fk R$ write $tHs$ [respectively $tVs$] to mean that there is a ``horizontal'' [respectively
``vertical''] directed edge
$[t, s]$ in $\cl G_h(\fk R)$ [respectively $\cl G_v(\fk R)$].
Define homomorphisms $T_1, T_2 : \bb Z \fk R \to \bb Z \fk R$ by
\begin{equation*}
  T_1t = \sum_{tHs}s, \qquad
  T_2t = \sum_{tVs}s.
\end{equation*}
It follows from (\ref{e2}),(\ref{e3}) that 
\begin{align*}
  T_1t &= \left(\sum_{b(s)=b'(t)}s\right) - t^h, \\
  T_2t &= \left(\sum_{a(s)=a'(t)}s\right) - t^v.
\end{align*}
Consider the homomorphism
\begin{align*}
  \left(\begin{matrix}T_1-I\\T_2-I\end{matrix}\right) : 
   \quad \bb Z \fk R & \to \bb Z \fk R  \oplus \bb Z \fk R,  \\
   t & \mapsto (T_1t-t)\oplus (T_2t-t).
\end{align*}
Define $\e : \bb Z \cl E \to \bb Z \cl E^+$ by 
\begin{equation*}
\e(x) = 
\begin{cases}
    \, x  & \text{if $x\in \cl E^+$},\\
    -\ovl x &  \text{if $x \in \overline{\cl E^+}$.}
\end{cases}
\end{equation*}
The boundary map  $\partial : \bb Z \fk R^+ \to \bb Z \cl E^+$ is defined by
$$\partial t = \e(a(t)+b'(t)-a'(t)-b(t))$$
and since $X$ is 2-dimensional, $H_2(\Gamma, \bb Z)=\ker \partial$.
Define a homomorphism 
$$\vf_2: \bb Z \fk R^+ \to \bb Z \fk R$$
 by
\begin{equation*}
  \vf_2t=t-t^v-t^h+t^{vh}.
\end{equation*}
The rest of this section is devoted to proving the following result, which is a more precise version
of Theorem \ref{main0}.
\begin{theorem}\label{main}
The homomorphism $\vf_2$ restricts to an isomorphism from $H_2(\Gamma, \bb Z)$ onto
  $\ker \left(\begin{smallmatrix}T_1-I\\T_2-I\end{smallmatrix}\right)$.
\end{theorem}
Define a homomorphism $\vf_1: \bb Z\cl E \to \bb Z \fk R \oplus \bb Z \fk R$ by
\begin{align*}
  \vf_1(a) &= 0\oplus\left(\sum_{a(s)=\ovl a}s-\sum_{a(s)=a}s \right), \qquad  \text{if $a\in \fk A$}, \\
  \vf_1(b) &= \left(\sum_{b(s)=b}s-\sum_{b(s)=\ovl b}s\right)\oplus 0, \qquad\,\,  \text{if $b\in \fk B$}.
\end{align*}
Note that if $x\in \cl E$ then $\vf_1(\ovl x)=-\vf_1(x)$ and so $\vf_1(\e(x))=\vf_1(x)$.

\begin{lemma}  The homomorphisms $\vf_1, \vf_2$ are injective and the following diagram commutes:
  \begin{equation}\label{cd}
\begin{CD}
\bb Z\cl E^+                   @<\partial<<   \bb Z \fk R^+\\
@V\vf_1VV                                              @VV\vf_2V \\
\bb Z \fk R \oplus \bb Z \fk R    @<<\left(\begin{smallmatrix}T_1-I\\T_2-I\end{smallmatrix}\right)<    \bb Z \fk R
\end{CD}
\end{equation}
\end{lemma}

\begin{proof}
Let $t\in \fk R$. Then
\begin{align*}
  (T_1-I)t &= \left(\sum_{b(s)=b'(t)}s\right) - t^h -t, \\
  (T_1-I)t^v &= \left(\sum_{b(s)=\ovl {b'(t)}}s\right) - t^{vh} -t^v, \\
  (T_1-I)t^h &= \left(\sum_{b(s)=b(t)}s\right) - t -t^h, \\
  (T_1-I)t^{vh} &= \left(\sum_{b(s)=\ovl{b(t)}}s\right) - t^v -t^{vh}. 
\end{align*}
Therefore
\begin{equation*}
\begin{split}
 (T_1-I)\circ\vf_2(t)& =  (T_1-I)(t-t^v-t^h+t^{vh})\\
& = \left(\sum_{b(s)=b'(t)}s-\sum_{b(s)=\ovl {b'(t)}}s\right)
-\left(\sum_{b(s)=b(t)}s-\sum_{b(s)=\ovl {b(t)}}s\right).
\end{split}
\end{equation*}
By definition of $\vf_1$, this implies that
$$
\vf_1(b'(t)-b(t))=(T_1-I)\vf_2(t)\oplus 0.
$$
Similarly
$$
\vf_1(a(t)-a'(t))=0\oplus (T_2-I)\vf_2(t).
$$
Therefore
\begin{equation*}
\begin{split}
\left(\begin{smallmatrix}T_1-I\\T_2-I\end{smallmatrix}\right)\circ \vf_2 (t)& =
\vf_1(b'(t)-b(t)+a(t)-a'(t))\\
& = \vf_1\circ \e(b'(t)-b(t)+a(t)-a'(t))\\
& = \vf_1\circ \partial (t).
\end{split}
\end{equation*}
This shows that (\ref{cd}) commutes.

It is obvious that $\vf_2$ is injective. To verify that $\vf_1$ is injective,
 define $\psi :  \bb Z \fk R \oplus \bb Z \fk R\to \bb Z\cl E^+$ by
  $\psi(s,t)=\e(b(s)-a(t))$.
Then $\psi\circ\vf_1(x)$ is a nonzero multiple of $x$, for all $x\in \cl E$. 
It follows that $\psi\circ\vf_1 : \bb Z\cl E^+ \to \bb Z\cl E^+$ is injective and therefore so is
$\vf_1$.
\end{proof}

\begin{lemma}\label{part}
  The homomorphism $\vf_2$ restricts to an isomorphism from $H_2(\Gamma, \bb Z)$ onto
  $\vf_2(\bb Z \fk R^+) \cap \ker \left(\begin{smallmatrix}T_1-I\\T_2-I\end{smallmatrix}\right)$.
\end{lemma}

\begin{proof}
  Let $\vf_2(\beta)\in \ker \left(\begin{smallmatrix}T_1-I\\T_2-I\end{smallmatrix}\right)$,
  where $\beta\in \bb Z \fk R^+$.  It follows from (\ref{cd}) that $\vf_1\circ\partial(\beta)=0$.
  But $\vf_1$ is injective, so $\partial \beta=0$ i.e. $\beta\in H_2(\Gamma, \bb Z)$. 
  
  Conversely, if $\beta\in H_2(\Gamma, \bb Z)$ then 
  $\left(\begin{smallmatrix}T_1-I\\T_2-I\end{smallmatrix}\right)\circ \vf_2 (\beta)=0$
by (\ref{cd}), so $\vf_2(\beta)\in \ker \left(\begin{smallmatrix}T_1-I\\T_2-I\end{smallmatrix}\right)$.
Since $\vf_2$ is injective, the conclusion follows.
\end{proof}

The next result, combined with Lemma \ref{part}, completes the proof of Theorem \ref{main}. 

\begin{lemma}\label{connected} There is an inclusion $\ker \left(\begin{smallmatrix}T_1-I\\T_2-I\end{smallmatrix}\right)\subset\vf_2(\bb Z \fk R^+)$.
\end{lemma}

\begin{proof}
Let $\alpha=\sum_{t\in \fk R} \lambda(t)t \in \ker \left(\begin{smallmatrix}T_1-I\\T_2-I\end{smallmatrix}\right)$.
We show that $\alpha \in \vf_2(\bb Z \fk R^+)$.
If $s\in \fk R$ then the coefficient of $s$ in the sum representing $(T_1-I)\alpha$ is
$$
\left(\sum_{\substack{t\in \fk R, t\ne s^h\\ b'(t)=b(s)}}\lambda(t)\right) - \lambda(s)
=\left(\sum_{\substack{t\in \fk R \\ b'(t)=b(s)}}\lambda(t)\right) - \lambda(s)-\lambda(s^h).
$$
This coefficient is zero, since $\alpha\in\ker (T_1-I)$. Therefore
\begin{equation}\label{l}
  \lambda(s)+\lambda(s^h)=\sum_{\substack{t\in \fk R \\ b'(t)=b(s)}}\lambda(t).
\end{equation}
The right hand side of equation (\ref{l}) depends only on $b(s)$, so for any $b\in \fk B$ we define
\begin{equation*}
\mu(b)=\sum_{\substack{t\in \fk R\\b'(t)=b}}\lambda(t).
\end{equation*}
Thus (\ref{l}) may be rewritten as
\begin{equation}\label{insert}
\lambda(s)+\lambda(s^h)=\mu(b(s)).
\end{equation}
It follows from (\ref{l}) and (\ref{flip}) that
\begin{equation}\label{induction}
\mu(b(s))=\mu(b(s^h))=\mu(b'(s)).
\end{equation}
\refstepcounter{picture}
\begin{figure}[htbp]
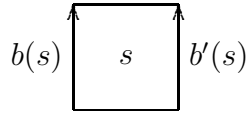

\hfil
\centerline{
\beginpicture
\setcoordinatesystem units <0.7cm, 0.7cm>  
\setplotarea  x from -6 to 6,  y from -1 to 1
\putrule from -1 -1 to 1 -1
\putrule from -1  1 to 1  1
\putrule from -1 -1 to -1 1
\putrule from  1  -1 to 1  1
\put {$b(s)$} [r]     at   -1.2 0
\put {$b'(s)$} [l]     at    1.2 0
\put {$s$}     at   0  0
\arrow <6pt> [.3,.67] from    -1  0.8   to  -1 1 
\arrow <6pt> [.3,.67] from   1  0.8    to   1  1
\endpicture
}
\hfil
\caption{$\mu(b(s))=\mu(b'(s))$}
\end{figure}

Fix an element $b_0\n \in \fk B$, and let $\cl C$ be the connected component of the graph $\cl G_h(\fk B)$ containing $b_0$. Then $\cl C$ is a connected graph with vertex set $\cl C^0 \subset \fk B$ and edge set $\cl C^1 \subset \fk R$. The graph $\cl C$ has a natural orientation $\cl C^+ = \cl C^1 \cap (\fk R^+\cup \fk R^v)$ and it is clear that  $\cl C^1= \cl C^+ \cup \{t^h : t\in \cl C^+\}$. Each vertex of $\cl C$ has degree at least three, since the same is true of the tree $\cl T_1$. Therefore the number of vertices of $\cl C$ is less than the number of geometric edges i.e. $|\cl C^0| < |\cl C^+|$.

If $b\in \cl C^0$ then there is a path in $\cl C^0$ from $b_0$ to $b$.
It follows by induction from (\ref{induction}) that $\mu(b_0)=\mu(b)$.
Thus
\begin{equation*}
\mu(b_0) =\sum_{\substack{t\in \fk R \\ b'(t)=b}}\lambda(t)=\sum_{\substack{t\in \cl C^1 \\ b'(t)=b}}\lambda(t).
\end{equation*}
Therefore 
\begin{equation*}
\begin{split}
|\cl C^0|\mu(b_0) & = \sum_{b\in \cl C^0}\sum_{\substack{t\in \cl C^1 \\ b'(t)=b}}\lambda(t)= \sum_{t\in \cl C^1}\lambda(t)\\
& = \sum_{t\in \cl C^+}(\lambda(t)+\lambda(t^h))= \sum_{t\in \cl C^+}\mu(b(t))\\
& = \sum_{t\in \cl C^+}\mu(b_0) = |\cl C^+|\mu(b_0).
\end{split}
\end{equation*}
Since $|\cl C^0| < |\cl C^+|$, it follows that $\mu(b_0)=0$ for all $b_0\in \fk B$. In other words, by (\ref{insert}),
\begin{equation}\label{h}
\lambda(s)=-\lambda(s^h)
\end{equation}
for all $s\in \fk R$.
A similar argument, using $\alpha\in\ker (T_2-I)$ and interchanging the roles of horizontal and vertical reflections, shows that
\begin{equation}\label{v}
\lambda(s)=-\lambda(s^v)
\end{equation}
for all $s\in \fk R$.
Combining (\ref{h}) and (\ref{v}) gives
\begin{equation}\label{vh}
\lambda(s)=\lambda(s^{vh})
\end{equation}
for all $s\in \fk R$.
Finally,
\begin{equation*}
\begin{split}
\alpha & =\sum_{t\in \fk R^+} \left(\lambda(s)s +\lambda(s^v)s^v + \lambda(s^h)s^h + \lambda(s^{vh})s^{vh}\right)\\
& = \sum_{t\in \fk R^+} \lambda(s)\left(s-s^v-s^h+s^{vh}\right)\\
&=\sum_{t\in \fk R^+} \lambda(s)\vf_2(s)  \in \vf_2(\bb Z \fk R^+).
\end{split}
\end{equation*}
\end{proof}

\section{K-theory of the boundary $C^*$-algebra}\label{five}

The (maximal) boundary $\partial\Delta$ of $\Delta$ is defined in \cite{kr}. It is homeomorphic to $\partial \cl T_1 \times \partial \cl T_2$, where $\partial \cl T_j$ is the totally disconnected space of ends of the tree $\cl T_j$.
The group $\Gamma$ acts on $\partial\Delta$ and hence on $C_{\bb C}(\partial\Delta)$ via $g\mapsto \alpha_{g}$, where 
$\alpha_{g} f(\omega)=f(g^{-1}\omega)$, for $f\in C_{\bb C}(\partial\Delta)$, $g\in\Gamma$.
The full crossed product $C^*$-algebra $\cl A(\Gamma,\partial\Delta)=C_{\bb C}(\partial\Delta) \rtimes \Gamma$ is the completion of the algebraic crossed product in an appropriate norm. We present examples where the rank of the analytic $K$-group $K_0(\cl A(\Gamma,\partial\Delta))$ is determined by Theorem \ref{main}.

\subsection{One vertex complexes} The case where the quotient VH-T complex $X$ has one vertex was studied in \cite{kr}. The group $\Gamma$ acts freely and transitively on the vertices of $\Delta$ and $\cl A(\Gamma,\partial\Delta)$ is isomorphic to a rank-2 Cuntz-Krieger algebra, as described in \cite{rs1, rs2}. The proof of this fact given in \cite[Theorem 5.1]{kr}. It follows from \cite{rs1} that $\cl A(\Gamma,\partial\Delta)$ is classified by its K-theory. 
By the proofs of \cite[Proposition 4.13]{rs2} and \cite[Lemma 4.3, Theorem 5.3]{kr}, we have 
$$K_0(\cl A(\Gamma,\partial\Delta))=K_1(\cl A(\Gamma,\partial\Delta))$$
and 
$$\rank(K_0(\cl A(\Gamma,\partial\Delta)))= 2\cdot\dim\ker\left(\begin{smallmatrix}T_1-I\\T_2-I\end{smallmatrix}\right).$$ 
Together with Theorem \ref{main}, this proves
\begin{equation}\label{principal2}
\rank K_0(\cl A(\Gamma,\partial\Delta))= 2\cdot\rank H_2(\Gamma, \bb Z).
\end{equation}
This verifies a conjecture in \cite{kr}.

\subsection{Irreducible lattices in $\pgl_2(\bb Q_p)\times\pgl_2(\bb Q_\ell)$}

If $p, \ell$ are prime then the group $\pgl_2(\bb Q_p)\times\pgl_2(\bb Q_\ell)$ acts on the $\Delta=\cl T_{p+1} \times \cl T_{\ell+1}$ and on its boundary $\partial\Delta$, which can be identified with a 
direct product of projective lines $\bb P_1(\bb Q_p) \times \bb P_1(\bb Q_\ell)$.
Let $\Gamma$ be a torsion free irreducible lattice in $\pgl_2(\bb Q_p)\times\pgl_2(\bb Q_\ell)$.
Then $\Gamma$ acts freely on $\Delta$ and $\cl A(\Gamma,\partial\Delta)$ is a rank-2
Cuntz-Krieger algebra in the sense of \cite{rs1}. The irreducibility condition (H2) of \cite {rs1} follows from Lemma \ref{nasc}. The proofs of the remaining conditions of \cite{rs1} are exactly the same as in \cite[Lemma 4.1]{kr}. It follows that (\ref{principal2}) is also true in this case.
Since $\Gamma$ is irreducible, the normal subgroup theorem \cite[IV, Theorem (4.9)]{mar} implies that $H_1(\Gamma, \bb Z)=\Gamma/[\Gamma, \Gamma]$ is finite.
Equation (\ref{principal2}) can therefore be written
\begin{equation}\label{fourteen}
\chi(\Gamma) = 1+ \frac{1}{2}\rank K_0(\cl A(\Gamma,\partial\Delta)).
\end{equation}
On the other hand, one easily calculates $\chi(\Gamma)=\frac{(p-1)(\ell-1)}{4}|X^0|$, where $|X^0|$
is the number of vertices of $X$. Therefore the rank of $K_0(\cl A(\Gamma,\partial\Delta))$ 
can be expressed explicitly in terms of $p,\ell$ and $|X^0|$. 

Explicit examples are studied in \cite[Section 3]{moz}.  If $p,l \equiv 1 \pmod 4$ are two distinct primes, Mozes constructs an irreducible lattice $\Gamma_{p,\ell}$ in $PGL_2(\bb Q_p) \times
PGL_2(\bb Q_l)$ which acts freely and transitively on the vertex set of $\Delta$. 
Here is how $\Gamma_{p,l}$ is constructed. Let $\bb H(\bb Z)=\{a=a_0+a_1i+a_2j+a_3k ; a_j\in \bb Z\}$, the ring of integer quaternions, let $i_p$ be a square root of $-1$ in $\bb Q_p$ and define 
$$\psi : \bb H(\bb Z) \to PGL_2(\bb Q_p) \times PGL_2(\bb Q_\ell)$$
by
$$\psi(a)=\left(
\begin{bmatrix}
a_0+a_1i_p & a_2+a_3i_p \\
-a_2+a_3i_p  & a_0-a_1i_p \\
\end{bmatrix},
\begin{bmatrix}
a_0+a_1i_\ell & a_2+a_3i_\ell \\
-a_2+a_3i_\ell  & a_0-a_1i_\ell \\
\end{bmatrix}
\right).
$$

Let $\tilde\Gamma_{p,\ell}=\{a=a_0+a_1i+a_2j+a_3k\in \bb H(\bb Z) ; a_0\equiv 1 \pmod 2,
a_j\equiv 0 \pmod 2, j=1,2,3, |a|^2=p^rl^s\}$.
Then $\Gamma_{p,\ell}=\psi(\tilde\Gamma_{p,\ell}$). The fact that $\Gamma_{p,\ell}$ is irreducible
follows easily from \cite[Corollary 2.3]{rr}, where it is observed that the only nontrivial direct product subgroup of $\Gamma_{p,\ell}$ is $\bb Z \times \bb Z = \bb Z^2$.

Since $|X^0|=1$, it follows from (\ref{fourteen}) that
the rank of $K_0(\cl A(\Gamma_{p,\ell},\partial\Delta))$ is
$\frac{(p-1)(\ell-1)}{2}-2$. This proves an experimental observation of \cite[Example 6.2]{kr}.
The construction of Mozes has been generalized in \cite[Chapter~3]{rat}
to all pairs $(p,l)$ of distinct odd primes and the same conclusion applies.

\end{document}